\documentclass{amsart}
\usepackage{amssymb,amscd}
\usepackage[all,2cell,ps]{xy}

\newcommand{\R}{{\mathbb{R}}}
\newcommand{\Q}{{\mathbb{Q}}}
\newcommand{\N}{{\mathbb{N}}}

\newcommand{\I}{\mathcal{I}}
\newcommand{\F}{{\mathcal{F}}}
\newcommand{\px}{\mathcal{P}(X)}
\newcommand{\pn}{\mathcal{P}(\N)}
\newcommand{\epx}{\text{End}(\mathcal{P}(X))}

\newcommand{\edge}[1]{\ar@{-}[#1]}

\newtheorem{theorem}{Theorem}[section]

\newtheorem{lemma}[theorem]{Lemma}

\theoremstyle{definition}
\newtheorem{definition}[theorem]{Definition}
\newtheorem{question}[theorem]{Question}

\theoremstyle{remark}

\numberwithin{equation}{section}

\begin{document}

\title[Kuratowski's 14-set theorem]{Variations on Kuratowski's 14-set theorem}
\author{David Sherman}
\address{Department of Mathematics\\ University of Illinois\\ Urbana, IL 61801}
\email{dasherma@math.uiuc.edu}

\begin{abstract}
Kuratowski's 14-set theorem says that in a topological space, 14 is the maximum possible number of distinct sets which can be generated from a fixed set by taking closures and complements.  In this article we consider the analogous questions for any possible subcollection of the operations \{closure, complement, interior, intersection, union\}, and any number of initially given sets.  We use the algebraic ``topological calculus" to full advantage.
\end{abstract}

\maketitle

\section{Introduction}

The following well-known result, from Kuratowski's 1920 dissertation, is known as the \textit{14-set theorem}.

\begin{theorem} \label{T:k} \cite{K}
Let $E \subset X$ be subset of a topological space.  The number of distinct sets which can be obtained from $E$ by successively taking closures and complements (in any order) is at most 14.  Moreover, 14 can be achieved if $X$ contains a subset homeomorphic to the Euclidean line.
\end{theorem}


The main goal of this article is to see what happens when ``closure" and ``complement" are replaced or supplemented with other basic topological operations.

\begin{question} \label{q}
Let $\I$ be a subcollection of
$$\{\text{closure, interior, complement, intersection, union}\}.$$
What is the maximum number of distinct sets which can be generated from a single set in a topological space by successive applications of members of $\I$?
\end{question}

Apparently Question \ref{q} will require us to answer $2^5 = 32$ different questions... well, not really.  Many of these are redundant, either because different choices of $\I$ allow us to perform the same operations, or because different choices of $\I$ raise algebraically isomorphic questions.  (This is explained at the end of Section \ref{S:answers}.)  Theorem \ref{T:k} answers at least one case, and certainly many others are trivial.  Finally, an example of Kuratowski shows that if we allow all five operations, we may obtain infinitely many sets.  After a full reckoning, there remain only two new questions to be answered.

\begin{question} ${}$ \label{Q}

\begin{enumerate}
\item What is the maximum number of distinct sets that can be generated from a fixed set in a topological space by successively taking closures, interiors, and intersections (in any order)?
\item Same question, but with closures, interiors, intersections, and unions.
\end{enumerate}
\end{question}

With a little additional work, we will finally answer

\begin{question} \label{q'}
Same question as Question \ref{q}, with $n \ge 2$ sets initially given. 
\end{question}

\bigskip

Our approach to this topic, like Kuratowski's, is almost entirely \textit{algebraic}.  The basic language is the ``topological calculus" which was developed by the Polish school during the first half of the twentieth century.  Prominent figures such as Birkhoff, Stone, Halmos, and Tarski kept this program dynamic, intertwining topology with the related fields of set theory, logic, and lattice theory.  (And the incorporation of Hilbert spaces opened up new realms of noncommutative analysis, with von Neumann at the center.)  In the present article, the topological calculus that we need is the type of universal algebra known as ``closure algebra," which reduces the questions above to the calculation of certain algebras generated by a specific partially ordered set.  So while our subject is apparently point-set topology, points and sets play a very minor role!

Theorem \ref{T:k} is actually easy to prove, and almost certainly has acquired some cachet from the unusual presence of the number 14.  Many people have extended Kuratowski's result by abstracting the algebraic content or isolating the specific conditions which allow a topological space and subset to generate 14 sets.  Some of their papers are mentioned in Section \ref{S:remarks} and the References.  The investigation closest to our own is due to Zarycki \cite{Z}, who proved some results analogous to Theorem \ref{T:k} by replacing ``closure" with other unary topological operations.  The ideas of allowing the binary operations of intersection and/or union (but not all three Boolean operations), and subsequently permitting $n \geq 2$ initial sets, seem to be new.  It turns out that our approach, unlike Zarycki's, produces some other ``unusual" numbers.

The paper is intended for the nonspecialist in universal algebra - indeed, it was written by one.  Thus we define even basic terms, and do not always give the most general formulations.  It is hoped that many readers will find the methods at least as interesting as the answers.

\section{Monoids, posets, and the proof of Theorem \ref{T:k}} \label{S:pf}

We start with the basics.  Let $X$ be an arbitrary set, and let $\px$ be the set of subsets of $X$.  To endow $X$ with a \textit{topology} means to choose a distinguished subset of $\px$, called the \textit{open sets}, which is closed under arbitrary unions and finite intersections, and contains both $X$ and the empty set.  The complement of an open set is a \textit{closed set}.  The (topological) \textit{closure} of $E \in \px$ is the smallest closed set containing $E$; the \textit{interior} of $E$ is the largest open set contained in $E$.  Therefore the three functions ``closure of," ``interior of," and ``complement of" can naturally be viewed as operations on $\px$.  We will denote them by $k$, $i$, $c$, respectively, and write them to the left of the set, as is usual for operators (or English sentences).  We also denote the collection of maps $\px \to \px$ as $\epx$.  Thus $kiE$ should be read as ``the closure of the interior of $E$."  The reader should be aware that some authors place topological operations to the right of the set (with an opposite rule for composition), and the letters $k$ and $c$ are sometimes switched.  With regard to the latter, our choice was made with Kuratowski closure operators in mind. 

\begin{definition} \label{D:kc} \cite{K}
A \textbf{Kuratowski closure operator} on a set $X$ is a map $k \in \epx$ which satisfies, for any $E,F \in \px$,
\begin{enumerate}
\item $k \emptyset = \emptyset$;
\item $kkE = kE$;
\item $kE \supseteq E$;
\item $kE \cup kF = k(E \cup F)$.
\end{enumerate}
\end{definition}

A Kuratowski closure operator $k$ on $X$ is exactly the topological closure operator for the topology on $X$ whose open sets are $\{ckE \mid E \subseteq X\}$.  So the choice of $k$ is equivalent to the choice of topology on $X$, while $c$ is independent of topology.  We write $I \in \epx$ for the identity map and record the following identities and their consequences:
\begin{equation} \label{E:ident}
k^2 = k, \quad c^2 = I, \quad i = ckc \quad \Rightarrow \quad i^2 = i, \quad ic = ck, \quad kc = ci.
\end{equation}

Now we recall some definitions from algebra.  A \textit{monoid} is a set with an associative binary operation and a unit.  (So a monoid is a ``group without inverses.")  A \textit{partial order} on a set is a reflexive antisymmetic transitive relation.  A standard example is $\leq$ on $\R$, but it is not necessary that any two elements be comparable: $\px$ is a partially ordered set - a \textit{poset} - with partial order given by inclusion.  One notices that $k$ and $i$ preserve the ordering, while $c$ reverses it.  (This means, for example, that $E \supseteq F \Rightarrow kE \supseteq kF$ - use Definition \ref{D:kc}(4).)  Now the class of functions from any set into a poset can also be made into a poset, where one function dominates another iff this is true pointwise.  Thus we have an induced partial order on $\epx$: for $\varphi, \psi \in \epx$,
$$\varphi \geq \psi \quad \iff \varphi(E) \supseteq \psi(E), \quad \forall E \in \px.$$
Then item (3) of Definition \ref{D:kc} can be rewritten as $k \ge I$, and apparently $i \le I$.  Note that order is preserved by an arbitrary right-composition: 
$$\varphi \geq \psi \Rightarrow \varphi \sigma \geq \psi \sigma, \qquad \varphi, \psi, \sigma \in \epx.$$
Order is also preserved by left-composition with $k$ or $i$, but reversed by left-composition with $c$.  

``Ordered monoid" sounds frighteningly abstract, but we will only be concerned here with subsets of $\epx$, with the binary operation of composition and the ordering $\leq$ as above.  The advantage in the situation at hand is that we may invoke a familiar friend from group theory (or universal algebra, to those in the know): \textit{presentations}.  This just means that we will describe sets of operations in terms of generators and relations, as demonstrated in

\begin{lemma} \label{L:ki} ${}$

\begin{enumerate}
\item Let $k,i \in \epx$ be the closure and interior operators of a topological space.  Then the cardinality of the monoid generated by $k$ and $i$ is at most $7$.
\item For a subset of a topological space, the number of distinct sets which can be obtained by successively taking closures and interiors (in any order) is at most $7$.
\end{enumerate}
\end{lemma}

\begin{proof} Composing $k \geq I$ on the left and right with $i$ gives $iki \geq i$.  Composing $i \leq I$ on the left and right with $k$ gives $kik \leq k$.  We use both of these to calculate
$$(i)k \leq (iki)k = i(kik) \leq i(k) \Rightarrow ik = ikik.$$
$$k(i) \leq k(iki) = (kik)i \leq (k)i \Rightarrow ki = kiki.$$
Since $k^2 = k$ and $i^2 = i$, the monoid generated by $k$ and $i$ contains exactly
$$\{I, i, ik, iki, k, ki, kik\}$$
(which may not all be distinct).  This proves the first part, and the second part is a direct consequence.
\end{proof}

We will use parentheses for ``the monoid generated by," so the second sentence of Lemma \ref{L:ki}(1) can be rewritten as $|(k,i)| \le 7$.

\begin{proof}[Proof of Theorem \ref{T:k}.] It follows from \eqref{E:ident} that any word in $k,i,c$ can be reduced to a form in which $c$ appears either as the leftmost element only, or not at all.  So by the previous lemma
\begin{equation} \label{E:14list}
(k,c) = (k,i,c) = \{I, i, ik, iki, k, ki, kik, c, ci, cik, ciki, ck, cki, ckik\}.
\end{equation}
Thus 14 is an upper bound.  To conclude the proof, it suffices to exhibit a so-called \textit{(Kuratowski) $14$-set}: a subset of a topological space for which all of these operations produce distinct sets.  One example is $S = \{0\} \cup (1,2) \cup (2,3) \cup [\Q \cap (4,5)] \subset \R.$
\end{proof}

\bigskip

We now investigate the order structure of $(k,i)$ a little further.  Using our basic rules for order we find that
$$i \leq I \leq k; \qquad i \leq iki \leq [\text{either of }ki, ik] \leq kik \leq k.$$
By considering the set $S$ we see that these (and the consequences from transitivity) are the only order relations, at least when $X$ contains a copy of $\R$.  Then the order structure of $(k,i)$ is depicted in Figure 1.  Here a segment from $\varphi$ up to $\psi$ means that $\varphi < \psi$ and there is no $\sigma$ satisfying $\varphi < \sigma < \psi$.  (We write $\varphi < \psi$ for $\varphi \leq \psi$  and $\varphi \ne \psi$.)

\begin{figure}[b]
\centerline{
\xymatrix{
            & & k\edge{ddll}\edge{dr}\\
            & &                     & kik\edge{dl}\edge{dr}\\
I\edge{ddrr}& & ki\edge{dr}         &                     & ik\edge{dl}\\
            & &                     & iki\edge{dl}\\
            & & i \\
}
}
\caption{The order structure of $(k,i)$.}
\end{figure}
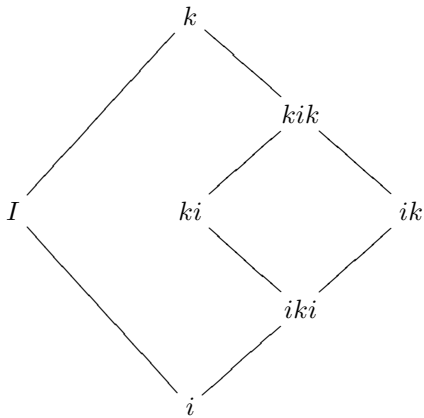

Since left composition with $c$ reverses order, and generically there are no relations between the first seven and last seven elements of \eqref{E:14list}, the order structure of $(k,i,c)$ consists of two disjoint copies of Figure 1.  It was first drawn by Kuratowski \cite{K}; in essence all of the arguments in this section go back to his dissertation.

\section{Lattices and freeness} \label{S:iu}

We have already mentioned that $\px$ is a poset.  Now we want to take advantage of its much richer structure as a \textit{Boolean lattice}.

Recall that a \textit{lattice} is a poset in which any two elements have a least upper bound and a greatest lower bound.  We write these binary operations as $\vee$ and $\wedge$, respectively, and refer to them as \textit{join} and \textit{meet}.  A lattice $L$ is \textit{distributive} if it satisfies the equality
$$x \wedge (y \vee z) = (x \wedge y) \vee (x \wedge z), \qquad x,y,z \in L$$
(or equivalently, the same equality with $\vee$ and $\wedge$ everywhere interchanged).  A lattice $L$ is \textit{complemented} if it contains a least element 0, a greatest element 1, and
$$ \forall a \in L, \; \exists b \in L \qquad a \vee b = 1, \quad a \wedge b = 0.$$
The element $b$ is called a \textit{complement} of $a$.  A Boolean lattice, finally, is a complemented distributive lattice.

It should be clear that $\px$ is a Boolean lattice with the operations of union, intersection, and set complementation.  As the class of all functions from $\px$ into the Boolean lattice $\px$, $\epx$ also inherits structure as a Boolean lattice, with pointwise operations.  It will be sufficient here to note the definitions of meet and join: for $\varphi, \psi \in \epx$, $E \in \px$, we define
$$(\varphi \wedge \psi) E = (\varphi E) \cap (\psi E), \qquad (\varphi \vee \psi) E = (\varphi E) \cup (\psi E).$$

We will be interested in Boolean lattices which are equipped with a closure operator $k$ satisfying the lattice version of Definition \ref{D:kc}.  Such an object is called a \textit{closure algebra}, and by \cite[Theorem 2.4]{MT} it is always isomorphic to a Boolean sublattice of some $\px$, where $X$ is a topological space and $k$ the associated closure operator.  (In broad generality a \textit{(universal) algebra} is a set equipped with operations \cite{Gua}.)  We define $i = ckc$ and note in particular Definition \ref{D:kc}(4), which says that $k$ distributes over $\vee$ (and $i$ distributes over $\wedge$).  So a closure algebra is a set equipped with the three unary operations $k,i,c$ and two binary operations $\wedge,\vee$, satisfying a certain small list of relations.  Any meaningful composition of these operations, for example $(E,F) \mapsto kE \wedge i(kcF \vee E)$, will also be called an operation.  We do not introduce the partial order relation formally but instead view $\varphi \leq \psi$ as shorthand for $\varphi \wedge \psi = \varphi$.

In the situation at hand we can specialize further to \textit{singly-generated} closure algebras.  That is, we suppose that every element can be obtained from a certain fixed generator (corresponding to the initial set of Question \ref{Q}) by some unary operation.  Of particular interest is the unique (up to isomorphism) \textit{free} singly-generated closure algebra \cite[Theorem 5.1]{MT}, referred to here as $\F$.  Freeness means, roughly, that the set of relations is minimal.  In other words, if two unary operations agree on the generator of $\F$, they agree on every element of every closure algebra.  This has the convenient consequence that inside $\F$, we need not distinguish between elements and unary operations.

Although not essential to the sequel, we pause to mention some striking facts about $\F$ and the logical structure of the theory of closure algebra.  They are taken from a highly-recommended 1944 article of McKinsey and Tarski \cite[Theorem 5.10, Theorem 5.17, Appendix IV]{MT}.
\begin{itemize}
\item $\F$ is isomorphic to a sub-closure algebra of the closure algebra of the Euclidean line, so that a certain subset of the line distinguishes all unequal unary operations of closure algebra.
\item There is an algorithm (guaranteed to terminate in finitely many steps) for deciding whether two operations of closure algebra are equal.
\item When two operations of closure algebra are equal, there is a formal proof.
\end{itemize}
Logicians also know closure algebra as a \textit{meta}mathematical object, since it can be used as a framework for intuitionistic logic \cite[Chapter III]{RS}.

\bigskip

From all this it follows that Question \ref{Q} is equivalent to
\begin{question} ${}$ \label{Q'}
\begin{enumerate}
\item What is the cardinality of $(k,i,\wedge)_\F$, the subalgebra of $\F$ generated by $\{k,i,\wedge\}$? 
\item Same question for $(k,i,\wedge,\vee)_\F$.
\end{enumerate}
\end{question}

\section{Answers to Questions \ref{q} and \ref{Q}} \label{S:answers}

We first determine the two algebras of Question \ref{Q'}.  Our strategy is simply to begin with the seven-element poset $(k,i)$ and add all meaningful compositions.

\bigskip

\textbf{1. Closures, interiors, intersections:}

We need to add all irredundant meets to Figure 1.  First we add $ki \wedge ik$.  Now notice that the meet of any two elements, both different from $I$, is already in our poset.  It suffices to add the meet of each element with $I$, and since $k \wedge I = I$ and $i \wedge I = i$, this gives five more elements.  The resulting structure is the \textit{meet semi-lattice} generated by the poset $(k,i)$, since it has $\wedge$ but not $\vee$.  In fact it is the same as the poset of all hereditary subsets of Figure 1, ordered by inclusion.  (A subset $H$ of a poset is \textit{hereditary} if $x \in H, \: x \geq y \Rightarrow y \in H$.)

A diagram of this 13-element poset is given in Figure 2.  By construction it is closed under $\wedge$, and since $i$ distributes across $\wedge$ it is closed under $i$.  Perhaps surprisingly, it is also closed under $k$.  To prove this, we need to show that for each element $\varphi$ in Figure 2, $k\varphi$ already appears in Figure 2.  

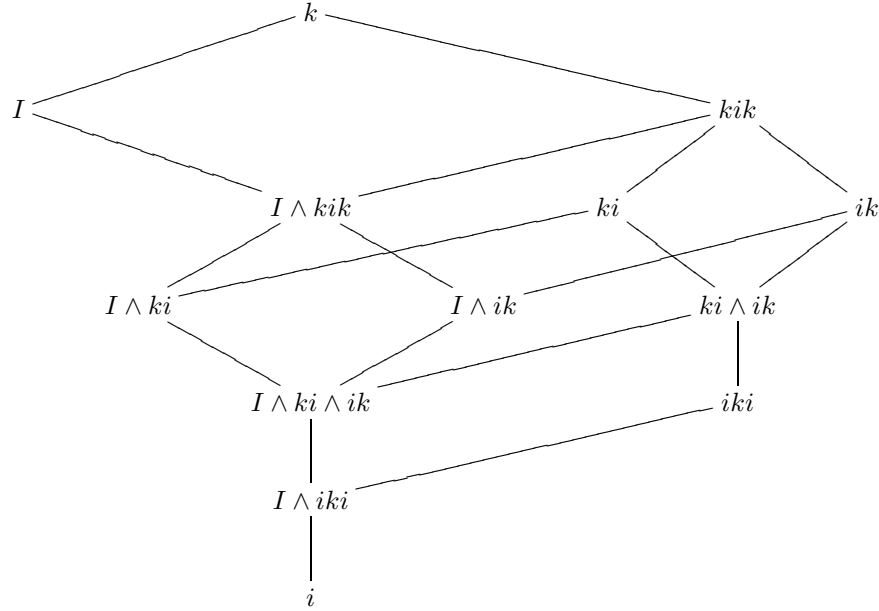
\begin{figure}
\centerline{
\xymatrix{
& & k \edge{dll} \edge{drrr} \\
I \edge{drr} &&&&& kik \edge{dlll} \edge{dl} \edge{dr} \\
&& I \wedge kik \edge{dl} \edge{dr}&&ki \edge{dr} \edge{dlll} && ik \edge{dlll} \edge{dl}\\                    
&I \wedge ki \edge{dr}&&I \wedge ik \edge{dl} && ki \wedge ik \edge{dlll} \edge{d} \\
&&I\wedge ki \wedge ik \edge{d} &&& iki \edge{dlll} \\
&&I \wedge iki \edge{d} \\
&&i \\ 
}
}
\caption{The order structure of $(k,i,\wedge)_\F$.}
\end{figure}

This is clear for the seven elements from Figure 1.  For the remaining elements, we start with an easy observation.  Since $k$ preserves order, for any $E,F \in \px$ we have
$$k(E \cap F) \subseteq kE, \; k(E \cap F) \subseteq kF \Rightarrow k(E \cap F) \subseteq kE \cap kF.$$
This means that 
\begin{equation} \label{E:kw}
k(\varphi \wedge \psi) \leq k\varphi \wedge k\psi, \qquad \varphi, \psi \in \epx.
\end{equation}
Now let $\sigma$ be any of $ki \wedge ik, I \wedge iki, I \wedge ki \wedge ik, I \wedge ki$.  Applying \eqref{E:kw} to $\sigma$ and reducing gives $k \sigma \leq ki$.  But $\sigma \geq i$, so $k \sigma \geq ki$.  We conclude that $k \sigma = ki$.

It is left to consider the two elements $(I \wedge ik)$ and $(I \wedge kik)$.  Applying \eqref{E:kw} shows that the closure of each is $\le kik$.  We claim that $k(I \wedge ik) = kik$, whence the larger element $k(I \wedge kik)$ is $kik$ as well.  We calculate as follows:
\begin{align*}
ik &= ik \wedge k(I) \\
&= ik \wedge k[(I \wedge ik) \vee (I \wedge cik)] \\
&= ik \wedge [k(I \wedge ik) \vee k(I \wedge cik)] \\
&= [ik \wedge k(I \wedge ik)] \vee [ik \wedge k(I \wedge cik)].
\end{align*}
Inspecting the last term,
$$ik \wedge k(I \wedge cik) \le ik \wedge k(cik) = ik \wedge cik = 0.$$
(Here $0 \in \epx$ is the map which sends every set to the empty set.)  We may therefore omit this term from the previous equation, which gives
$$ik = ik \wedge k(I \wedge ik) \Rightarrow ik \leq k(I \wedge ik) \Rightarrow kik \leq k(I \wedge ik).$$
Since the opposite inequality was already established, the claim is proved.

It follows that $(k,i,\wedge)_\F$ has at most thirteen elements, and it can be checked that each of these operations produces a distinct set when applied to the set
\begin{align} \label{E:14set}
T &= \left[ \bigcup_{n \in \N} \left\{ \frac{1}{n} \right\} \right] \cup \left[ [2,4] - \bigcup_{n \in \N} \left\{3 + \frac{1}{n} \right\} \right] \\
\notag &\quad \cup \left[ (5,7] \cap \left( \Q \cup \bigcup_{n \in \N} \left(6 + \frac{1}{2n\pi}, 6 + \frac{1}{(2n-1)\pi} \right] \right) \right].
\end{align}
So the answer to Question \ref{Q}(1) is \textit{thirteen}.  The set $T$ also shows that there are no order relations in $(k,i,\wedge)_\F$ other than those of Figure 2, which is important for the next question.

\bigskip

\textbf{2. Closures, interiors, intersections, unions:}

Our first task here is to add all irredundant joins to Figure 2.  We will consider expressions of the form
\begin{equation} \label{E:meetirr}
x_1 \vee x_2 \dots \vee x_n, \qquad x_j \in (k,i,\wedge)_\F.
\end{equation}
Since $\F$ is distributive, the resulting set will be closed under joins and meets: it is in fact the distributive lattice generated by Figure 1.

Let us add in the two elements $ki \vee ik$ and $(I \wedge ki) \vee (I \wedge ik)$, and partition our poset into four classes:
\begin{enumerate}
\item $i,k$;
\item $I$;
\item $iki, ki \wedge ik, ik, ki, ki \vee ik, kik$;
\item $I \wedge iki, I \wedge ki \wedge ik,  I \wedge ik, I \wedge ki, (I \wedge ki) \vee (I \wedge ik), I \wedge kik$.
\end{enumerate}
It may help to notice that the third class is the right-hand five of Figure 2, plus $ki \vee ik$, while the fourth class is the middle five of Figure 2, plus $(I \wedge ki) \vee (I \wedge ik)$.

Each class above is already a sublattice of $\F$, so an irredundant join of the form \eqref{E:meetirr} can contain at most one $x_j$ from each class.  Elements in the first class occur in no irredundant joins.  $I$ cannot be involved in an irredundant join except with elements of the third class, which produces six more elements.  It is left to consider joins of the third and fourth classes.  Using the distributive law, this turns up 14 more elements:
\begin{align*}
&(I \wedge kik) \vee ki \vee ik, (I \wedge kik) \vee ki, (I \wedge kik) \vee ik, \\
&\qquad (I \wedge kik) \vee (ki \wedge ik), (I \wedge kik) \vee iki, \\
&(I \wedge ki) \vee (I \wedge ik) \vee (ki \wedge ik), (I \wedge ki) \vee (I \wedge ik) \vee iki, \\
&(I \wedge ki) \vee ik, (I \wedge ki) \vee (ki \wedge ik), (I \wedge ki) \vee iki, \\
&(I \wedge ik) \vee ki, (I \wedge ik) \vee (ki \wedge ik), (I \wedge ik) \vee iki, \\
&(I \wedge ki \wedge ik) \vee iki.
\end{align*}
We conclude that the distributive lattice generated by the poset $(k,i)$ has at most 35 elements, each of which is a join of elements from Figure 2.  Since $k$ distributes across joins, this set is closed under $k$.  The roles of $k$ and $i$ are dual (see below for explanation) in the distributive lattice generated by $(k,i)$, so it is also closed under $i$.  Finally, one can check that the 35 operations are distinguished by the set $T$ from \eqref{E:14set}.  Therefore the answer to Question \ref{Q}(2) is \textit{thirty-five}.

\bigskip

\begin{table}[b]
\begin{center}
\begin{tabular}{|l|r|r|r|r|}
\hline
Operations & $\{I\}$ & $\{\wedge\}$ & $\{\vee\}$ & $\{\wedge, \vee\}$ \\ \hline
$\{I\}$ & 1 & 1 & 1 & 1 \\ \hline
$\{i\}$ & 2 & 2 & 2 & 2 \\ \hline
$\{k\}$ & 2 & 2 & 2 & 2 \\ \hline 
$\{c\}$ & 2 & 4 & 4 & 4 \\ \hline 
$\{i,k\}$ & 7 & 13 & 13 & 35 \\ \hline
$\{i,c\} = \{k,c\} = \{i,k,c\}$ & 14 & $\infty$ & $\infty$ & $\infty$ \\ \hline
\end{tabular} \caption{Solution to Question \ref{q}.  Each number is the cardinality of the subalgebra of $\F$ generated by the operations in its row and column.}
\end{center}
\end{table}

A complete answer to Question \ref{q} is given in Table 1.  All of the numbers $\leq 4$ in Table 1 are trivial to verify, and some of the repetition is due to the fact that in the presence of $c$, the inclusion of $k$ or $i$ (resp. $\wedge$ or $\vee$) is equivalent to the inclusion of $k$ and $i$ (resp. $\wedge$ and $\vee$).  Other repetition is due to duality, which we now describe.

The dual of a poset is the same underlying set, with the ordering reversed.  (So the diagram is turned upside-down.)  The Boolean lattice $\px$ is isomorphic with its own dual, via the complementation map $c$.  So any operation $\varphi$ on $\px$ has a dual operation, $\bar{\varphi} = c\varphi c_n$, where by $c_n$ we mean the application of $c$ to each of the $n$ arguments of $\varphi$.  One pictures $\bar{\varphi}$ as a vertical reflection of $\varphi$; $k$ and $i$ are dual, as are $\wedge$ and $\vee$.  The duality mapping $\varphi \mapsto \bar{\varphi}$ commutes with composition and so is an automorphism of $\epx$.  It restricts to an automorphism of $\F$ (which is closed under duality, as is any subalgebra which contains $c$).  In particular it preserves the cardinality of subsets.  This explains, for example, why $|(k,i,\wedge)_\F| = |(k,i,\vee)_\F|$.

Finally we give a simplified version of Kuratowski's example \cite{K} showing that $\F = (k,c,\wedge)_\F$ is infinite.  Define a closure operator $k$ on $\pn$ by
\begin{equation} \label{E:pn}
k(A) =
\begin{cases}
\varnothing, & A = \varnothing; \\
[\min A, \infty) = \{\min A, 1 + \min A, \dots\}, & \text{otherwise;}
\end{cases}
\end{equation}
for any $A \in \pn$.  It is easy to see that $k$ satisfies Definition \ref{D:kc} and so determines a topology on $\N$.  Let $\varphi = I \wedge [k(k \wedge c)]$, and let $E \subset \N$ be the even numbers.  The reader can easily check that $\varphi^j(E) = E \cap [2j+2, \infty)$.  Since these sets are all distinct, the operations $\varphi^j$ are distinct, and $|\F| = \infty.$ 

\section{Answers to Question \ref{q'}}

Now we alter our hypotheses by supposing that $n \ge 2$ sets are initially given.  The theorems of McKinsey and Tarski apply to this case as well, so that we may consider the free closure algebra generated by $n$ elements.  We will denote it as $\F_n$, and the generators as $\{F_n\}$.  Remarkably, $\F_n$ also embeds in the closure algebra of the Euclidean line.

At first glance Question \ref{q'} may seem intractable or at least extremely tedious.  In the presence of $\wedge$ or $\vee$ the cardinalities grow at least exponentially: for example, the two sets $T \times \R, \R \times T \subset \R^2$ obviously generate more than $13^2 = 169$ distinct subsets under $k,i,\wedge$.  It turns out, however, that all of the hard work is done, and the situation stabilizes nicely for $n \ge 2$.  A complete solution to Question \ref{q'} is given in Table 2; below we explain the key points.

\begin{table}[b]
\begin{center}
\begin{tabular}{|l|r|r|r|r|}
\hline
Operations & $\{I\}$ & $\{\wedge\}$ & $\{\vee\}$ & $\{\wedge, \vee\}$ \\ \hline
$\{I\}$ & $n$ & $2^n - 1$ & $2^n - 1$ & $D_n$ \\ \hline
$\{i\}$ & $2n$ & $3^n - 1$ & $\infty$ & $\infty$ \\ \hline
$\{k\}$ & $2n$ & $\infty$ & $3^n - 1$ & $\infty$ \\ \hline 
$\{c\}$ & $2n$ & $2^{2^n}$ & $2^{2^n}$ & $2^{2^n}$ \\ \hline 
$\{i,k\}$ & $7n$ & $\infty$ & $\infty$ & $\infty$ \\ \hline
$\{i,c\} = \{k,c\} = \{i,k,c\}$ & $14n$ & $\infty$ & $\infty$ & $\infty$ \\ \hline
\end{tabular} \caption{Solution to Question \ref{q'}.  Each number is the cardinality of the subalgebra of $\F_n$ generated by the operations in its row and column.}
\end{center}
\end{table}

We first claim that $|(k,\wedge)_{\F_2}| = \infty$, proved in almost the same way as $|\F| = \infty$.  Let $k$ be the closure operator on $\pn$ of equation \eqref{E:pn}, and let $E = E_0$ and $O$ be the even and odd positive integers, respectively.  For $j \ge 1$, define inductively elements of the closure algebra generated by $E$ and $O$ by
$$E_j = E_{j-1} \cap k(kE_{j-1} \cap O).$$
Then the $E_j = E \cap [2j+2, \infty)$ are all distinct, establishing the claim.  By duality and inclusions, this justifies every occurrence of $\infty$ in Table 2.

The first column of Table 2 consists of algebras with unary operations only, so the results of Table 1 can be applied to one generator at a time.  For the last three entries in the fourth row of Table 2, the algebra under consideration is the free Boolean algebra with $n$ generators.  It is well-known to have $2^{2^n}$ elements \cite[Theorem II.2.2(iii)]{G}.

The elements of $(\wedge)_{\F_n}$ have the form $F_{j_1} \wedge F_{j_2} \dots \wedge F_{j_k}$, where $1 \le k \le n$ and the $j_k \in \{1,2,\dots n\}$ are distinct.  Apparently they are in 1-1 correspondence with the $2^n - 1$ nonempty subsets of the $n$ generators.  The situation for $(i, \wedge)_{\F_n}$ is similar; now any element is a meet in which for each $F_j$, one of three possibilities holds: $F_j$ is present, $iF_j$ is present, or both are absent.  (Recall that $i$ distributes across $\wedge$.)  Since we do not admit the empty meet, this allows $3^n - 1$ possibilities.  The other occurrences of $2^n - 1$ and $3^n - 1$ in Table 2 follow from duality.

Finally, $(\wedge, \vee)_{\F_n}$ is the free distributive lattice on $n$ generators; determining its cardinality $D_n$ is sometimes known as \textit{Dedekind's problem}.  No explicit formula for $D_n$ is known, but there has been much work on its asymptotics.  For the curious reader, we recall Kleitman's result \cite{Kl}: $\log_2 D_n \sim C(n, \lfloor\frac{n}{2}\rfloor)$, where $C(n,k)$ denotes the binomial coefficient and $\lfloor x \rfloor$ the greatest integer $\le x$.

Notice that each of the finite formulas of Table 2 extends to the case $n=1$.

\section{Remarks and related questions} \label{S:remarks}

Our approach suggests many other questions related to Kuratowski's 14-set theorem.  For example, one can consider different subalgebras of $\F$, in which the generating operations include other topological operations built out of $\{k,\wedge,c\}$.  This need not be too abstract -- for example, the operation ``boundary of" is $k \wedge kc$.  Some cases are discussed by Zarycki \cite{Z}.

Other basic topological operations arise which do not belong to closure algebra.  We briefly mention the operation ``accumulation points of," which sends $E \in \px$ to its \textit{derived set} $DE$.  The operator $D$ is not expressible in terms of $k,\wedge,c$, as consideration of the set $\{0\} \cup \{1/n\}$ shows.  This observation goes back to Kuratowski \cite{K}.

It is also possible to generalize Definition \ref{D:kc} by relaxing the postulates and/or replacing $\px$ with non-Boolean lattices.  Among the sizable research in this direction, we point out \cite{H}, which deals specifically with Kuratowski's 14-set theorem.

The present article is about the algebraic content of Kuratowski's 14-set theorem, but there are also investigations into its topological content.  We have seen that a topological space $X$ contains 14-sets when it is sufficiently rich, in particular when it contains a copy of the Euclidean line.  One can ask for characterizations of topological spaces in which certain limitations occur, and a classification along these lines was given by Aull \cite{A}.  For a specific subset $E \subseteq X$, Chapman \cite{C} gave an exhaustive description of all possible degeneracies of the poset in Figure 1, and Langford \cite{L} and Shum \cite{S} gave necessary and sufficient conditions for $E$ to be a 14-set.  It might be interesting (or it might not) to consider this kind of question in the setup of the questions considered above.

\bigskip

\textbf{Acknowledgment.} The author has benefited from several conversations with John D'Angelo.

\end{document}